\documentclass[12pt]{article}

\usepackage{geometry} 
\usepackage{amsthm,amssymb}
\usepackage{amsmath,amsfonts,latexsym}
\usepackage{latexsym}
\usepackage{float}
\usepackage{xcolor}
\usepackage{url}

\newtheorem{thm}{Theorem}[section]

\theoremstyle{definition}
\newtheorem{rem}[thm]{Remark}

\begin{document}

\title{New symmetric $2$-$(176,50,14)$ designs}
\author{
Dean Crnkovi\'c (\url{deanc@math.uniri.hr})\\
and\\
Andrea \v{S}vob (\url{asvob@math.uniri.hr})\\
Department of Mathematics \\
University of Rijeka \\
Radmile Matej\v ci\'c 2, 51000 Rijeka, Croatia
}

\maketitle

\begin{abstract}
In this paper we construct two new symmetric designs with parameters 2-(176,50,14) as designs invariant under certain subgroups of the full automorphism group of the Higman design. One is self-dual and has the full automorphism group of size 11520 and other is not self-dual and has the full automorphism group of size 2520.
\end{abstract}

\section{Introduction and preliminaries}

A $2$-$(v,k, \lambda)$ design is a finite incidence structure $(\mathcal{P}, \mathcal{B}, \mathcal{I})$, where $\mathcal{P}$ and $\mathcal{B}$ are disjoint sets and 
$\mathcal{I} \subseteq \mathcal{P} \times \mathcal{B}$, with the following properties:
\begin{enumerate}
  \item [1.] $|\mathcal{P}|=v$ and $1< k <v-1$,
  \item [2.] every element (block) of $\mathcal{B}$ is incident with exactly $k$ elements (points) of $\mathcal{P}$,
	\item [3.] every two distinct points in $\mathcal{P}$ are together incident with exactly $\lambda$ blocks of $\mathcal{B}$.
\end{enumerate}

In a $2$-$(v,k,\lambda )$ design every point is incident with exactly $\displaystyle r=\frac{\lambda (v-1)}{k-1}$ blocks, and $r$ is called the replication number of a design.
A 2-design is also called a block design. The number of blocks in a block design is denoted by $b$. If $v=b$, a design is called \textit{symmetric}. In a symmetric design, every two blocks meet in exactly $\lambda$ points. In case of symmetric designs, the notation $(v,k,\lambda)$ design is often used. An isomorphism from one block design to another is a bijective mapping of points to points and blocks to blocks which preserves incidence. An isomorphism from a block design ${\mathcal{D}}$ onto itself is called an automorphism of ${\mathcal{D}}$. The set of all automorphisms forms a group called the full automorphism group of ${\mathcal{D}}$ and is denoted by $Aut({\mathcal{D}})$. If $\mathcal{D}$ is a block design, the incidence structure $\mathcal{D}'$ having as points the blocks of $\mathcal{D}$, and having as blocks the points of $\mathcal{D}$, where a point and a block are incident in $\mathcal{D}'$ if and only if the corresponding block and a point of $\mathcal{D}$ are incident, is a block design called the {\it dual} of $\mathcal{D}$. The design is called self-dual if the design and its dual are isomorphic. For further reading on block designs we refer the reader to \cite{bjl}.

In \cite{higman} G. Higman has discovered the unique symmetric design with parameters 2-(176,50,14) having the Higman-Sims simple group as the full automorphism group. More detailed description of this design has been given in \cite{sims}, while details about strongly regular graph with parameters (176,49,12,14) connected with the Higman's design was done by Brouwer in \cite{brouwer}. The linear codes associated with Higman design are studied in \cite{par-ton}. The Higman symmetric 2-(176,50,14) design is self-dual.

In 1995, Z. Janko \cite{janko} constructed new symmetric design with parameters 2-(176,50,14) not isomorphic to the one constructed by Higman. The constructed design is self-dual and has the full automorphism group of order 1344.

In 1999, Z. Bo\v zikov \cite{zdravka} constructed another new symmetric design with parameters 2-(176,50,14) not isomorphic to the ones constructed before. The constructed design is not self-dual and has the full automorphism group of order 1000.

Up to our best knowledge, these four designs were up to isomorphism the only known designs with this parameters. In this paper we describe a construction which gives us three new symmetric 2-(176,50,14) designs, one self-dual design and one pair of dual designs.

We used Magma \cite{magma} for all computations involving groups and codes in this paper. The adjacency matrices of the three previously known designs and the two newly found designs, so up to isomorphism and duality all known symmetric 2-(176,50,14) designs, are available online at
\begin{verbatim}
 http://www.math.uniri.hr/~asvob/designs176.txt
\end{verbatim}

\section{Construction of new symmetric (176,50,14) designs}

The full automorphism group of the Higman design is the Higman-Sims simple group of order 44352000. In this section, we use a method for constructing symmetric designs based on orbit matrices that was developed in \cite{janko-orb}. We use the idea of taking some subgroups of the full automorphism group of the known design and then determine all symmetric designs 2-(176,50,14) on which the taken subgroup acts as an automorphism group in the same way as it acts on the original design.

Let $\mathcal{D}$ be a Higman design and let $G$ denote the Higman-Sims group, which is the full automorphism group of $\mathcal{D}$.

There is exactly one conjugacy classes of subgroups of order 11520 in the group $G=\mathrm{Aut}(\mathcal{D})$ of order 44352000, the representative of which will be denoted by $H_1$. It is isomorphic to the group $(E_{16}:A_6):Z_2$.
The subgroup $H_1$ acts in two orbits on the set of points and the set of blocks of $\mathcal{D}$, one of size 80 and the other of size 96, giving the orbit matrix $M_1$.

\begin{displaymath}  
\label{om1}
M_1=   \left(
\begin{tabular}{rr}
26  & 24 \\
20 & 30 \\
\end{tabular}  \right)
\end{displaymath}

The orbit matrix $M_1$  expands to the  (0,1)-incidence matrices
of two non-isomorphic symmetric designs: the design $\mathcal{D}$, and a new design, denoted by $\mathcal{D}_1$ with the full automorphism  group of order 11520, isomorphic to the group $H_1$.
The design $\mathcal{D}_1$ is self-dual.

There are exactly three conjugacy classes of subgroups of order 2520 in the group $G=\mathrm{Aut}(\mathcal{D})$ of order 44352000, the representative of which will be denoted by $H_2, H_3$ and $H_4$. 
The subgroups $H_2$ and $H_3$ act on the set of points (and the set of blocks) of $\mathcal{D}$ in three orbits, while $H_4$ acts in four orbits. The subgroups $H_2, H_3$ and $H_4$ are isomorphic to the group $A_7$. We will present a construction of a symmetric $2$-$(176,50,14)$ design from the group $H_2$. 

The subgroup $H_2$ acts in three orbits on the set of points of $\mathcal{D}$, of sizes 1, 70 and 105, and it acts in three orbits on the set of blocks of $\mathcal{D}$, of sizes 15, 35 and 126, giving the orbit matrix $M_2$.

\begin{displaymath}  
M_2=   \left(
\begin{tabular}{rrr}
1  & 28  &  21 \\
1  & 16  &  33 \\
0  & 20  &  30 \\
\end{tabular} \right)
\end{displaymath}

The orbit matrix $M_2$  expands to the  (0,1)-incidence matrices of two non-isomorphic symmetric designs:
the design $\mathcal{D}$, and a new design, denoted by $\mathcal{D}_2$ with  full automorphism group of order 2520, isomorphic to the group $H_2$. The design $\mathcal{D}_2$ is not self-dual,
we denote its dual design with $\mathcal{D}_3$. The design $\mathcal{D}_3$ can be constructed directly using the subgroup $H_3$.

The subgroup $H_4$ does not produce any new symmetric $2$-$(176,50,14)$ design. i.e., it produces only the Higman design.
	
\bigskip

The results presented above can be summarized as follows.


\begin{thm}
Up to isomorphism, there are exactly two symmetric $2$-$(176,50,14)$ designs invariant under the subgroup $H_1$ of order $11520$ of the Higman-Sims group acting with respect to the orbit matrix $M_1$, the Higman design and the design denoted by $\mathcal{D}_1$. The design $\mathcal{D}_1$ is self-dual, and $H_1$ is its full automorphism group.
Further, up to isomorphism there are exactly two symmetric $2$-$(176,50,14)$ designs invariant under the subgroup $H_2 \cong A_7$ of the Higman-Sims group acting with respect to the orbit matrix $M_2$, the Higman design and the design denoted by $\mathcal{D}_2$.  The full automorphism group of the design $\mathcal{D}_2$ is $H_2$. The design $\mathcal{D}_2$ is not self-dual, we denote its dual design with $\mathcal{D}_3$.
\end{thm}

\begin{rem}
There are at least seven symmetric $2$-$(176,50,14)$ designs, the Higman design, the design constructed by Janko, two mutually dual designs constructed by Bo\v zikov, and the designs $\mathcal{D}_1$, 
$\mathcal{D}_2$ and $\mathcal{D}_3$. 
\end{rem}

It is well known that designs with high degree of symmetry could produce codes with low dimension. If $p$ is a prime, $p$-rank of the incidence matrix of a $2$-$(v,k,\lambda)$ design can be smaller than 
$v-1$ only if $p$ divides the order of a design, i.e. if $p$ divides $k-\lambda$ in the case of a symmetric design (see \cite{hamada}). For a symmetric (176,50,14) design, the primes that could be of interest are 2 and 3. In Table \ref{results} we give the information about all known symmetric designs with parameters $2$-$(176,50,14)$.

\begin{table}[H]
\begin{center} \begin{footnotesize}
\begin{tabular}{|c|c|c|c|c|c|}
\hline
 $|Aut(D)|$& $Aut(D)$ & Self-dual & 2-rank & 3-rank &note\\
\hline
\hline
44352000 & $HS$  &yes & 22&50 & Higman, \cite{higman} \\
11520 & $(E_{16}:A_6):Z_2$   &yes & 46&66 & new\\
2520&$A_7$   &no & 26&51 & new \\
1344 &$(Z_4\times Z_4\times Z_4)\cdot F_{21}$     &yes & 48&50 & Janko, \cite{janko}\\
1000 &$(E_{25}: Z_{5}):Z_8$    &no & 62&70 & Bo\v zikov, \cite{zdravka}\\
\hline
\hline
\end{tabular}\end{footnotesize} 
\caption{\footnotesize Symmetric (176,50,14) designs, up to isomorphism and duality}\label{results}
\end{center}
\end{table}

\vspace*{0.5cm}
\begin{center}{\bf Acknowledgement}\end{center}
This work has been fully supported by {\rm C}roatian Science Foundation under the projects 6732 and 5713.


\begin{thebibliography}{30}

\bibitem{bjl}
T. Beth, D. Jungnickel, H. Lenz, Design Theory, 2nd Edition, Cambridge University Press, Cambridge, 1999.

\bibitem{magma}
W.~Bosma, J.~Cannon, Handbook of Magma Functions, Department of Mathematics, University of Sydney, 1994.
http://magma.maths.usyd.edu.au/magma.

\bibitem{zdravka}
Z. Bo\v zikov, A new symmetric design with parameters (176,50,14), J. Combin. Des. 8 (1999), 387--390.

\bibitem{brouwer}
A. E. Brouwer, Polarities of G. Higman's symmetric design and a strongly regular graph on 176 vertices, Aequationes Math. 25 (1982), 77--82.

\bibitem{hamada}
N. Hamada, On the $p$--rank of the incidence matrix of a balanced or partially
balanced incomplete block design and it applications to error correcting codes, Hiroshima Math. J. 3 (1973), 153--226.

\bibitem{higman}
G. Higman, On the simple group of D. G. Higman and C. C. Sims, Illinois J. Math. 13 (1969), 74--84.

\bibitem{janko-orb}
Z. Janko, Coset Enumeration in Groups and Constructions of Symmetric Designs, Combinatorics '90 (Gaeta, 1990), Ann. Discrete Math. 52 (1992), 275--277.

\bibitem{janko}
Z. Janko, On symmetric designs with parameters (176,50,14), J. Combin. Theory Ser. A 72 (1995), 310--314.

\bibitem{par-ton}
C. Parker, V. Tonchev, Linear Codes and Doubly Transitive Symmetric Designs, Linear Algebra Appl. 226--228 (1995), 237--246.

\bibitem{sims}
C. C. Sims, On the isomorphism of two groups of order 44 352000, Theory of Finite Groups (Symposium, Harvard Univ., Cambridge, 1968), Benjamin, New York, 1969, 101--108.

\end{thebibliography}
\end{document}